\documentclass[review]{elsarticle}
\usepackage{xcolor}
\usepackage{lineno,hyperref}
\modulolinenumbers[5]
\journal{Computers and Industrial Engineering}

\usepackage{amsmath}
\usepackage{amssymb}
 \usepackage{amsthm}
 \usepackage{amsfonts}
\usepackage{multirow}
\usepackage{eqnarray}
\usepackage{diagbox} 
\usepackage{eqnarray}
\newcommand{\rP}{\mathrm{P\space }} %probability, when spacing isn't a worry
\newcommand{\rE}{\mathrm{E\space }} %expectation, when spacing isn't a worry
\usepackage{float}
\newtheorem{theorem}{Theorem}

\newtheorem{definition}[theorem]{Definition}

\newtheorem{proposition}[theorem]{Proposition}

\begin{document}

\begin{frontmatter}

\title{Optimal maintenance schedule for a wind turbine with aging components}

\author[1]{Quanjiang Yu\corref{mycorrespondingauthor}} 
\address[1]{Department of Mathematical Sciences, Chalmers University of Technology and  University of Gothenburg, SE-42196 Gothenburg, Sweden}
\ead{yuqu@chalmers.se}
%\author[mysecondaryaddress]{Quanjiang Yu\corref{mycorrespondingauthor}}
\author[2]{Ola Carlson}
\address[2]{Department of Electrical Engineering, Chalmers University of Technology } %wait to be fill in
\author[1]{Serik Sagitov}
\cortext[mycorrespondingauthor]{Corresponding author}
\begin{abstract}
Wind power is one of the most important sources of renewable energy. A large part of the wind energy cost is due to the cost of maintaining the wind power equipment. To further reduce the maintenance cost, one can improve the design of the wind turbine components. One can also reduce the maintenance costs by optimal scheduling of the component replacements. The latter task is the main motivation for this paper. 

When a wind turbine component fails to function, it might need to be replaced under less than ideal circumstances. This is known as corrective maintenance.
To minimize the unnecessary costs, a more active maintenance policy based on the life expectancy of the key components is preferred. Optimal scheduling of preventive maintenance activities requires advanced mathematical modeling. 

In this paper, an optimization model is developed using the renewal-reward theorem. In the multi-component setting, our approach involves a new idea of virtual maintenance which allows us to treat each replacement event as a renewal event even if some components are not replaced by new ones.

The proposed optimization algorithm is applied to a four-component model of a wind turbine and the optimal maintenance plans are computed for various initial conditions. The modeling results showed clearly the benefit of PM planning compared to pure CM strategy (about $8.5\%$ lower maintenance cost). When we compare it with another state-of-art optimization model, it shows similar scheduling with a much faster CPU time. The comparison demonstrated that our model is both fast and accurate.
\end{abstract}


\end{frontmatter}

\section{Introduction}  
Wind power technology is one of the most efficient sources of the renewable energy available today.
A large part of the wind energy cost is due to the cost of maintaining the wind power equipment, especially for offshore wind farms.  
A corrective maintenance (CM) of a turbine component, performed after a break-down of the component, is usually more expensive than a preventive maintenance (PM) event, as some of the equipment is replaced in a planned manner. However, if PM activities are scheduled too frequently, the maintenance costs become unreasonably high, which entails the necessity of a maintenance schedule minimising the expected replacement, logistic, and downtime costs.

There is a broad body of literature devoted to various optimization models of maintenance scheduling. Here we name just few of relevant papers that influenced our own approach. In the article by \cite{lee2016new}, some general PM optimization models are presented.   
The effect of a PM action has been classified into three categories: failure rate reduction,
the decrease of the deterioration speed, and age reduction. 
The paper concludes that it can be profitable to perform PM.

 The article \cite{sarker2016minimizing} looks at opportunistic maintenance which is a special kind of preventive maintenance. When one component breaks down, the maintenance personal attending the broken component might as well maintain other aged components to save some logistic costs. This is extremely beneficial for offshore wind farms, due to the large set-up costs.

 In \cite{moghaddam2011sensitivity}, optimization models are developed to determine optimal PM schedules in repairable and maintainable systems. It was demonstrated that higher set-up costs  make advantageous simultaneous PM activities. However, the suggested models are nonlinear, which means they are computationally hard to solve.

The age dynamics of our model is based on the discrete Weibull distribution, see  \cite{guo2009reliability}. To account for the component deterioration due to aging, we assume that the PM cost increases as a linear function of the component's age at the replacement moment. Previously,  the PM replacement cost was usually treated as an age independent constant, see for example \cite{paper1}, \cite{santos2019modeling}.
Our approach focussing on age-dependence should be compared to that of \cite{chen2012bivariate}, where the maintenance cost is assumed to depend on the total damage.  Another relevant paper  \cite{liu2014value}, quantifies the maintenance cost of a component using the reliability distribution. 

Our definition of the objective function %$ f_{(s,\boldsymbol a)}(\boldsymbol x_s, \boldsymbol y_s)$ 
requires an application of the classical renewal-reward theorem, see for example   \cite{grimmett2020probability}. This approach is quite straightforward in the one component case as explained in  Section \ref{farm}. In the multiple component case, however, when  the true renewal events  (all components are replaced simultaneously) are rarely encountered, our approach requires some kind of pseudo-renewal events. To this end, in Section \ref{X} we introduce the key idea of  the virtual maintenance replacement cost $b_{(t,a)}$ for a generic component having age $a$ at time $t$.

Section \ref{CS} presents two sensitivity analyses and two case studies treating a four component model of the wind power turbine. In particular, it was demonstrated that under the additional assumptions  of \cite{paper1},  the new model produces similar results to those obtained in \cite{paper1} but at a higher computational speed. Some technical proofs of our claims are postponed until the end of the paper.

\section{A single component model}\label{farm}

\begin{figure}[h]
\centering
    \includegraphics[width=9 cm]{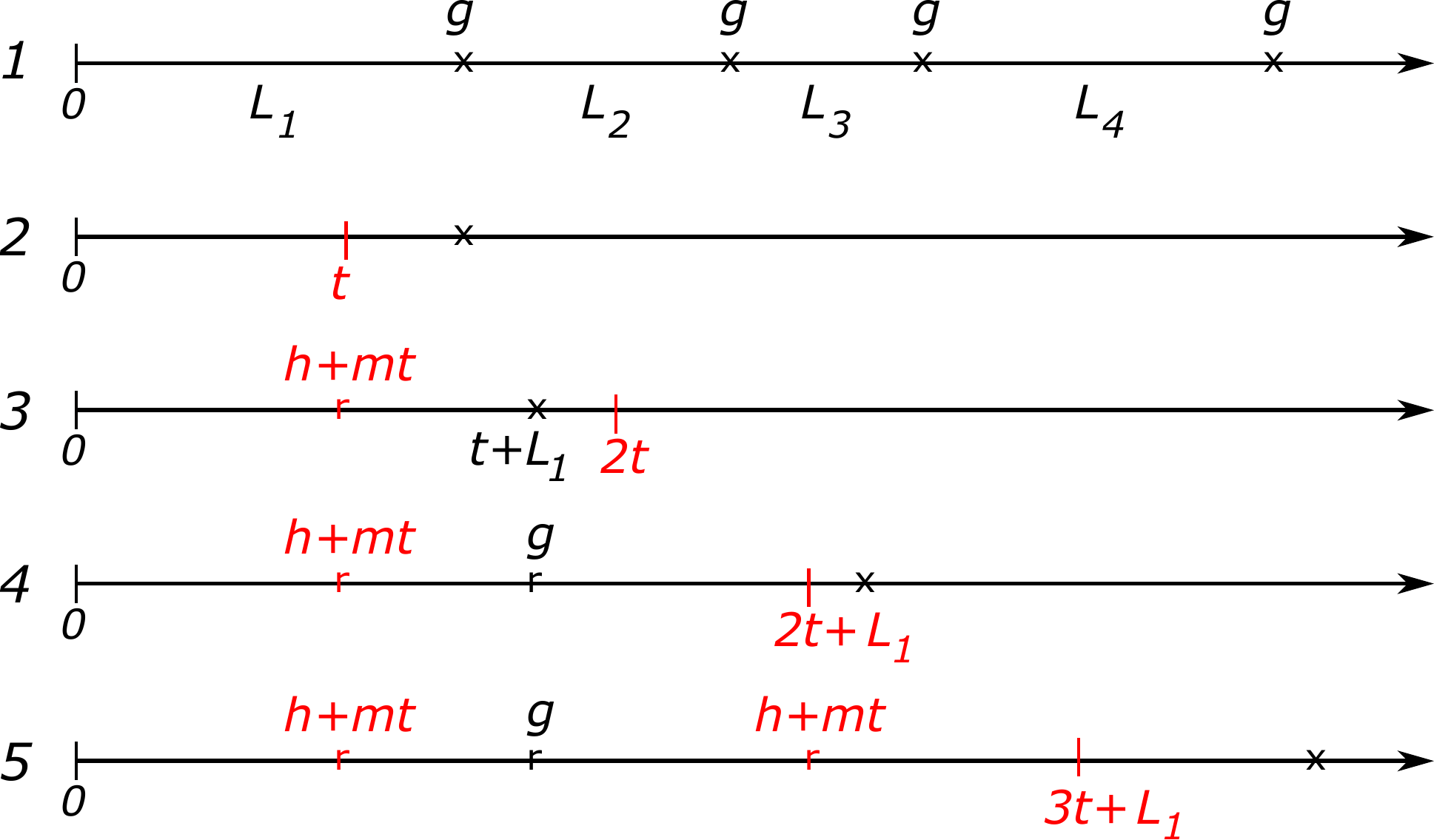}
    \caption{Renewal-reward model for computing $c$ for the one component model.}
    \label{Ser1}
\end{figure}

We start our exposition of the model by turning to the one component case. Without planned PM activities, the maintenance cost flow is described by line 1 of Figure \ref{Ser1}. 
Here, the consecutive failure times (depicted by crosses) form a renewal process with independent inter-arrival times $L_1, L_2,\ldots$ each having  the same Weibull distribution W$(\theta,\beta)$. After each failure, the broken component is replaced by a new one and the incurred replacement cost is $g$. According to the classical renewal-reward theorem the long term time-average reward (the maintenance cost in the current setting) is $c=\frac{g}{\rE(L)}$.
%\begin{description}
%%\item[ $a=$] the age of the gearbox at time $0$,
%%\item[ $e_t=$] the revenue generated by a turbine at time $t$,
%\item[ $g=$]  the corrective replacement cost,
%\item[ $h+mt=$] the age $t$ dependent preventive  replacement cost.
%\end{description}

The lines 2-5 on Figure \ref{Ser1} introduce the renewal-reward model used to compute the time-average maintenance cost  $q_t$ based on the PM planning strategy: plan the next PM replacement at time $t$ after each replacement. We use the red color to depict PM planning times, PM replacement times,  and PM costs; and the black color to depict the CM replacement times and CM costs.
Line 2 describes a scenario with the first failure occurring after the planned PM time and the first replacement is performed at time $t$, thereby the failure is avoided. The corresponding PM replacement cost is assumed to be a linear function $h+tm$ of the component's age at the time of a PM replacement. 

According to line 3, after the PM replacement was performed at the cost $h+tm$, the next failure occurred before the next PM planned time $2t$. As a consequence, the second replacement was performed in a CM regime incurring the cost $g$, see line 4. The line 4 says that following the second replacement, the next PM is planned at time $2t+L_1$, that is at time $t$ after the failure time  $t+L_1$. According to the placement of a cross on the line 4, we see on line 5 that the third replacement is performed in the PM regime and the next PM time is scheduled at time $3t+L_1$.

Recall that a random variable $L$ has a Weibull distribution W$(\theta,\beta)$ if
\[\rP(L>t)=e^{-\theta t^\beta},\quad t\ge 0,\]
where $\theta>0$ is a scale parameter and $\beta>1$ is a shape parameter of W$(\theta,\beta)$. The mean value of $L$ 
\[\mu=\theta^{-1/\beta}\Gamma(1+\tfrac{1}{\beta})\]
is computed using the gamma function.
The age dependence in the continous time setting is best viewed in terms of the hazard function $\theta \beta t^{\beta-1}$, 
%\[\lim_{\epsilon\to0}\tfrac{\rP(t<L<t+\epsilon|L>t)}{\epsilon}=\theta \beta t^{\beta-1}\]
telling how fast the failure rate increases with component's age $t$. In this paper, we use the discrete time version of the Weibull distribution 
\[\rP(L> t)=e^{-\theta t^\beta},\quad t=0,1,2,\ldots,\]
so that
\[\rP(L= t)=e^{-\theta (t-1)^\beta}-e^{-\theta t^\beta},\quad t=1,2,\ldots\]

Using these expressions and the renewal-reward theorem we arrive at the following explicit formula involving the key parameters of the one component model $(\theta,\beta,g,h,m)$. The proof of this result (as well as the forthcoming Proposition \ref{propc}) is given in the end of the paper.
\begin{proposition}\label{propq}
Think of an infinite planning horizon and a recurrent  strategy of planning the next PM at time $t$ after each replacement event. 
Then the  time-average maintenance cost is the following function 
\[q_t=\frac{(1-e^{-\theta t^\beta})g+e^{-\theta t^\beta}(h+mt)}{\sum_{k=1}^t (e^{-\theta (k-1)^\beta}-e^{-\theta k^\beta})k+e^{-\theta t^\beta}t}\]
of the planning time $t$.
\end{proposition}
(As a check one can send $t$ to infinity and observe that this results in $q_t\to\frac{g}{\rE(L)}$, the average maintenance cost in the absence of PM activities.)

Minimising the function $q_t$ over the possible planning times $t$, produces a constant
\begin{equation}
c=\min_{t\ge 1}\frac{(1-e^{-\theta t^\beta})g+e^{-\theta t^\beta}(h+mt)}{\sum_{k=1}^t (e^{-\theta (k-1)^\beta}-e^{-\theta k^\beta})k+e^{-\theta t^\beta}t},
 \label{c1}
\end{equation}
which we will treat as the long-term maintenance cost per unit of time for the component in question.
Notice that according to this formula,
$c$ is independent of the planning period $[s,T]$.

\section{The objective function in the one component case}\label{minP}
In this section we propose an optimization model for a PM scheduling during a discrete-time interval 
$$[s,s+1,\ldots,T-1,T]$$ 
 assuming that at time $s$ the component in use is of age $a$. If $a=0$, we say that at the beginning of the planning period the component was as good as new. This will allow us to find an optimal time $t_{(s,a)}$  for the next PM  by minimising the total maintenance cost during the whole planning period.

%The following definition brings a key ingredient of our optimization model.
\begin{definition}
 For the one component model with a planning period $[s,T]$, we call a {\rm PM plan} any vector
$$\boldsymbol x_s=(x_{s+1},\ldots,x_{T+1})$$ 
with binary components
%$x_t\in \{0,1\}$
satisfying a linear constraint
\begin{align}
\sum_{t=s+1}^{T+1}x_t=1,\quad x_{s+1}\in \{0,1\},\ldots, x_{T+1}\in \{0,1\}.\label{Con}
\end{align} 
\end{definition}
%We set $B_{T+1}=0$ and $\bar B_{T+1}=0$.

For the given planning period $[s,T]$, we define the total maintenance cost $Q_{(s,a)}(t,u)$ as a function of the planning time $t$ for the next PM and the failure time $u$ of the component in use.   For $t\in[s+1, T]$, put
\begin{align}
 Q_{(s,a)}(t,u)=\left\{
\begin{array}{ll}
g+(T-u)c , &  \text{if } u\leq t   \\
h+(t-s+a)m+(T-t)c,  &   \text{if } u> t 
\end{array}
\right.
\label{Qsa}
\end{align} 
or using indicator functions,
 \[
        Q_{(s,a)}(t,u)=[g+(T-u)c] 1_{\{u\leq t\}}  +  [h+(t-s+a)m+(T-t)c] 1_{\{u>t\}}.
\]
Put also
\[
        Q_{(s,a)}(T+1,u)=[g+(T-u)c] 1_{\{u\leq T\}}.
\]
This formula recognises two possible outcomes:
\begin{description}
\item[]  if $\{u\leq t\}$, then the breakdown happens before the planned PM time and the expected total maintenance cost is estimated to be $g+(T-u)c$, with $c$ given by \eqref{c1},
\item[] if  $\{u\geq t+1\}$, so that there is no breakdown before the planned PM time, then the expected total maintenance cost  is estimated to be  $$h+(t-s+a)m+(T-t)c.$$
\end{description}

If at the starting time $s$ of the planning period, the component in use has age $a\ge 1$, we will use a special notation $s+L_a$ for the first failure time, where
\[L_a\stackrel{d}{=}\{L-a|L>a\},\]  
is defined in terms of the full life length $L$  of a generic component. If $L$ has a discrete W$(\theta,\beta)$ distribution, then  
\[\rP(L_a> t)=\exp\left\{\theta\big(a^{\beta}-(a+t)^{\beta}\big)\right\},\quad t\ge0.\]

Putting $u=s+L_a$ into the formula for $Q_{(s,a)}(t,u)$, we arrive at  a random variable
\begin{align*}
    F_{(s,a)}(\boldsymbol x_s)=\sum_{t=s+1}^{T+1}Q_{(s,a)}(t,s+L_a)x_t
\end{align*}
that gives us the total maintenance cost of the PM plan $\boldsymbol x_s$. Averaging over $L_a$, we obtain the objective function
\[f_{(s,a)}(\boldsymbol x_s)=\rE(F_{(s,a)}(\boldsymbol x_s))\]
as the expected  maintenance cost of the PM plan $\boldsymbol x_s$. 
Now we are ready to define the optimal maintenance plan as the solution of the following optimization problem:
\begin{align*}
 \text{minimize}\qquad & f_{(s,a)}(\boldsymbol x_s) \\
\text{subject to}\qquad & \text{linear constraints }\eqref{Con}.
\end{align*}

Let $t_{(s,a)}$ be the PM time proposed by the solution of the minimisation problem above, and notice that 
\[t_{(0,0)}=\text{argmin } (q_t). \]
The following proposition states a consistency property for the set of optimal times $t_{(s,a)}$,  providing with an intuitive support for the suggested approach.
\begin{proposition}\label{propc}
Suppose for some positive $\delta$,
$$t_{(s,a)}> s+\delta.$$ 
If $t_{(s,a)}\le T$, then
\[t_{(s+\delta,a+\delta)}=t_{(s,a)}.\] 
\end{proposition}

\section{Multiple component model}\label{5.3}

In this section we expand our one component model to the $n\ge2$ component case. We now think of a wind turbine consisting of $n$ components with $j$-th component having a life length $L^j$ distributed according to a Weibull distribution W$(\theta^j,\beta^j)$, where parameters $(\theta^j,\beta^j)$ may differ for different $j=1,\ldots,n$. We will assume that components may require different replacement costs: 
\begin{description}
%\item[ $g_0=$] the shared logistic and down-time costs associated with a  CM activity,
\item[ $g_0=$] the down-time cost associated with a  CM activity,
\item[ $g^j=$] the component specific CM cost,
%\item[ $h_0=$] the fixed logistic cost plus the downtime cost during a PM activity,
\item[ $h_0=$] the downtime cost during a PM activity,
\item[ $h^j+tm^j=$] the component specific PM replacement cost for $j$-th component at age $t$.%I am not sure for this term whether we should add g_0 or not.
\end{description}

\begin{figure}[h]
\centering
    \includegraphics[width=9 cm]{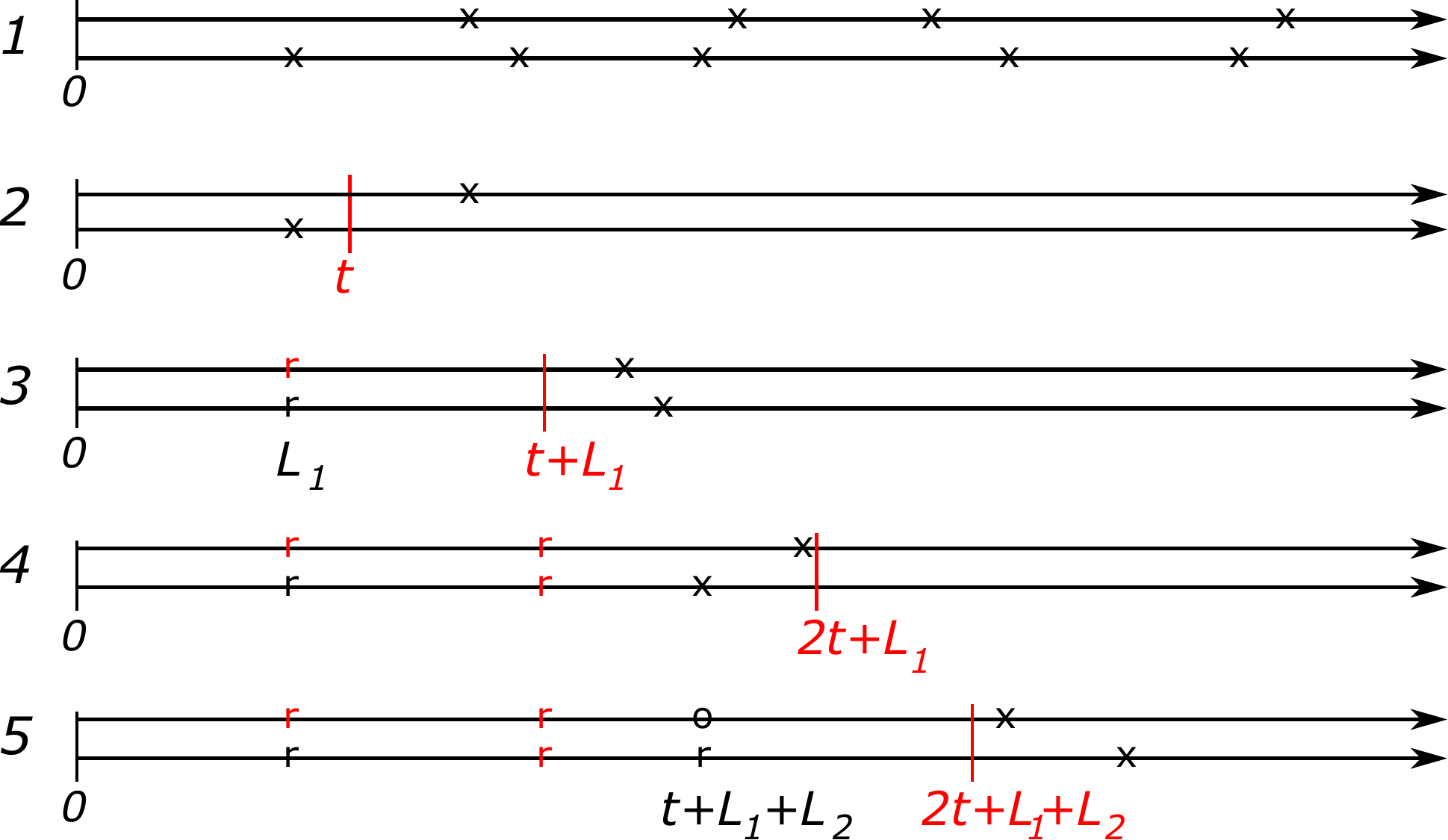}
    \caption{Renewal-reward model for computing $c$ for the two component model.}
    \label{Ser2}
\end{figure}

To be able to update the formula \eqref{c1} of the minimal time-average  maintenance cost, we would need to determine the times of total renewal of the system which are not readily available in the multiple component setting. For illustration turn to Figure \ref{Ser2} dealing with the case of $n=2$ components.
Even in the absence of PM activities, see line 1, it is clear that the classical renewal-reward theorem is not directly applicable. Our solution to this problem is to treat each replacement event as total renewal events, sometimes by performing opportunistic maintenance and in some cases, by increasing the cost function to reflect the future additional age-related replacement costs.
This idea is illustrated by lines 2-5 on Figure \ref{Ser2}.

Consider a strategy when the next PM activity is planned at time $t$ after each replacement. Line 2 depicts a case when the second component is broken first and before the time $t$ of the next PM. As shown by line 3, both components are replaced at the failure time $L_1$ and the next PM is planned at time $L_1+t$. In this way, the time $L_1$ can be viewed as the total renewal time of the two component system. The incurred replacement cost is
\[g_0+g^2+h^1+L_1m^1.\]

Since according to the line 3, both components break down after the planned time $t+L_1$ for the next PM, both components are replaced in the PM regime, so that the incurred replacement cost is
\[h_0+h^1+tm^1+h^2+tm^2.\]
Continuing in this way by replacing both components at each maintenance event, we arrive at a renewal-reward process described in the general setting as follows.

Assume that we start at time 0 with $n$ new components and denote by
\[L=\min(L^1,\ldots, L^n)\]
the time of the first failure. By independence, we have 
\[\rP(L>t)=\rP(L^1>t)\cdots\rP(L^n>t)=\exp\{-\theta^1 t^{\beta^1}-\ldots-\theta^n t^{\beta^n}\}.\]
Thus, if we plan our next PM at time $t$ after each maintenance replacement, the first renewal time is $X=X(t)$ with
\begin{equation}
X=L\wedge t=L\cdot 1_{\{L\leq t\}}+ t\cdot 1_{\{L>t\}} 
\label{Xt}
\end{equation}
and the corresponding reward value $\tilde R=\tilde R(t)$ is computed as
\[\tilde R=\Big(\sum_{j\in\gamma}g^j+\sum_{j\notin \gamma} (h^j+Lm^j)+g_0\Big) 1_{\{L\leq t\}}+ \Big(\sum_{j=1}^n (h^j+tm^j)+h_0\Big)1_{\{L> t\}},\]
where
\[\gamma=\{j: L^j=L\}\]
is the subset of components that would be replaced at the first failure. 
The renewal-reward theorem allows us to express the time-average maintenance cost $\frac{\rE(\tilde R)}{\rE(X)}$ as an explicit  function of the planning time $t$.

However, replacing all of the components at each maintenance event irrespectively of the ages of the components in use, is definitely a suboptimal strategy. If at a maintenance event some component $j$ is in a working condition and its current age  $a$ is rather small, then it might be more beneficial to let it continue working. For still being able to use the renewal argument in such a case , we introduce the idea of a virtual replacement. We replace the previous naive formula for reward $\tilde R$ by a more sophisticated one
\[R=\Big(\sum_{j\in\gamma}g^j+\sum_{j\notin \gamma} B^j_L+g_0\Big) 1_{\{L\leq t\}}+ \Big(\sum_{j=1}^n B^j_t+h_0\Big)1_{\{L> t\}},\]
where
\[B^j_a=(h^j+am^j)\wedge b^j_a\]
chooses the minimum between two age specific costs: the preventive replacement cost and the virtual replacement cost, see Section \ref{X}. 
The renewal-reward theorem implies that the time-average maintenance cost $\frac{\rE(R)}{\rE(X)}$ is computed as the following function of the planning time $t$
\[q_t=\frac{\rE_t(\sum_{j\in\gamma}g^j+\sum_{j\notin \gamma} B^j_L+g_0)+(\sum_{j=1}^n B^j_t+h_0) \rP(L> t)}
{\rE_t(L )+ t\rP(L> t)},\]
where $\rE_t(Z)$ stands for $\rE(Z\cdot 1_{\{L\leq t\}})$.
Now, after minimising $q_t$ over $t$ we define the desired constant
\begin{equation}
c=\min_{t\ge 1}\frac{\rE_t(\sum_{j\in\gamma}g^j+\sum_{j\notin \gamma} B^j_L+g_0) +(\sum_{j=1}^n B^j_t+h_0) \rP(L> t)}
{\rE_t(L )+ t\rP(L> t)}. \label{c}
\end{equation}

\section{Virtual replacement cost}\label{X}

Returning to the one component model and focussing on an arbitrary component $j$, observe that in this one component case the model parameters are $g=g_0+g^j$, $h=h_0+h^j$, 
$\theta=\theta^j$, $\beta=\beta^j$, and $m=m^j$. Using this set of parameters, we introduce the virtual replacement cost function $b_{(t,a)}^j=b_{(t,a)}$ as a function of the current time $t$ and the current age $a$ of the component in question. 

Since $f_{(s,a)}^*$ is the total cost of the optimal PM-plan in the one component setting. In terms of this cost function, 
the virtual replacement cost is defined by the difference
\begin{equation}
b_{(t,a)}=f_{(t,a)}^*-f_{(t,0)}^*.
\label{bta}
\end{equation}
This difference evaluates the extra maintenance cost over the time period $[t,T]$ due to the component's age $a$ at the starting time $t$ of the observation period.  The larger is $a$, the higher is the expected maintenance cost. The function 
$$b_a=b_{(0,a)}$$ 
will be treated as the long term virtual replacement cost of the component of age $a$.

\begin{proposition}\label{prop6}
If $t_{(0,0)}+s\leq T$, then $b_{(s,a)}=b_a$ and
\begin{equation}
t_{(s,a)}=t_{(0,0)}+s-a.
\label{tas}
\end{equation}

\end{proposition}

In the above described manner we can define $b_{(t,a)}^j=b_{(t,a)}$ and $b^j_a=b_a$ for all $j=1,\ldots, n$. 

\section{A PM plan and the objective function}\label{plan}
For the planning period $[s,T]$, we call a PM plan in the multicomponent setting any array $(\boldsymbol x_s,\boldsymbol y_s,z)$, where
$$(\boldsymbol x_s,\boldsymbol y_s)=\{x^{j}_{t},y_t: \ 1\le  j\le n,\ s+1\le t\le T\}$$  
and all components are binary $x^ {j}_{t},y_{t},z\in \{0,1\}$ subject to the following linear constraints
\begin{subequations}
\begin{align}
&y_t\geq x^{j}_{t},\quad t=s+1,\ldots T,\ j=1,\ldots,n,\label{tu}\\
&\sum_{t=s+1}^{T} y_t=1-z,\label{y}\\
&\sum_{i=1}^n x^j_t\geq y_t,\quad t=s+1,\ldots T.\label{xy}
\end{align} 
\label{const}
\end{subequations}
The equality $x^ {j}_{t}=1$ means that the plan $(\boldsymbol x_s,\boldsymbol y_s,z)$ suggests  a PM replacement of $j$-th component at time $t$. 
On the other hand, the equality
 $y_{t}=1$  means that according to the plan $(\boldsymbol x_s,\boldsymbol y_s,z)$  at least one of the components should be replaced at time $t$, this is guaranteed by constraint \eqref{xy} (constraint  \eqref{tu} allows for several components to be replaced at such a time $t$). The equality $z=1$ means that no PM is planned during the whole time period $[s+1,T]$. This is guaranteed by constraint \eqref{y}.

Given the current ages of $n$ components 
\[\boldsymbol a=(a^1,\ldots,a^n),\]
the first failure time is $s+L_{\boldsymbol a}$, where
$$L_{\boldsymbol a}=\min(L^1_{a^1},\ldots, L^n_{a^n}).$$
Putting
\[\gamma=\{j: L^j_{a^j}=L_{\boldsymbol a}\}\]
 we first mention a naive formula for the cost $ \tilde   F_{(s,\boldsymbol a)}$ assigned to a PM plan  $(\boldsymbol x_s,\boldsymbol y_s,z)$ assuming that at each maintenance event all $n$ components are replaced by the new ones. With
 \begin{align*}
 \tilde   C_{\boldsymbol a}&=g_0+(T-s-L_{\boldsymbol a}) c+  \sum_{j\in\gamma}g^j +\sum_{j\notin \gamma} (h^j+(L_{\boldsymbol a}+a^j)m^j),\\
\tilde   P_{\boldsymbol a}&=h_0+(T-t)c+\sum_{j=1}^{n}(h^j+(t-s+a^j)m^j),
\end{align*}
where $c$ is given by \eqref{c}, put
%\begin{align*}
%    \tilde F_{(s,\boldsymbol a)}(\boldsymbol y_s)
% =\sum_{t=s+1}^{T+1}y_t\Bigg[\Big(  \sum_{j\in\gamma}g^j& +g_0+(T-s-L_{\boldsymbol a}) c\\
%   &+\sum_{j\notin \gamma} (h^j+(L_{\boldsymbol a}+a^j)m^j)\Big)1_{\{s+L_{\boldsymbol a}\leq T\wedge t\}}\\
%   & \hspace{-2cm}+\Big(\sum_{j=1}^{n}(h^j+(t-s+a^j)m^j)+h_0+(T-t)c\Big)1_{\{s+L_{\boldsymbol a}> t\}}\Bigg],
%\end{align*}
\begin{align*}
   \tilde   F_{(s,\boldsymbol a)}(\boldsymbol y_s,z)
   =\sum_{t=s+1}^{T}\Big( &\tilde   C_{\boldsymbol a}1_{\{s+L_{\boldsymbol a}\leq t\}}+\tilde   P_{\boldsymbol a}1_{\{s+L_{\boldsymbol a}> t\}}\Big)y_t+\tilde   C_{\boldsymbol a}1_{\{s+L_{\boldsymbol a}\leq T\}}z.
   \end{align*}
Here the last term describes the option of planing no PM activity. 
This formula should be modified to incorporate the  virtual replacement costs:
%\[B^j_{(s,a)}=(h^j+am^j)\wedge b^j_{(s,a)}\]
%bringing
%\begin{align*}
%    F_{(s,\boldsymbol a)}(\boldsymbol y_s)
%   =\sum_{t=s+1}^{T+1}y_t\Bigg[\Big(  \sum_{j\in\gamma}g^j& +g_0+(T-s-L_{\boldsymbol a}) c\\
%   &+\sum_{j\notin \gamma} B^j_{(s+L_{\boldsymbol a},a^j+L_{\boldsymbol a})}\Big)1_{\{s+L_{\boldsymbol a}\leq T\wedge t\}}\\
%   & \hspace{-2cm}+\Big(\sum_{j=1}^{n}B^j_{(t,a^j+t-s)}+h_0+(T-t)c\Big)1_{\{s+L_{\boldsymbol a}> t\}}\Bigg].
%\end{align*}
\begin{align*}
    F_{(s,\boldsymbol a)}(\boldsymbol y_s,z)
    =\sum_{t=s+1}^{T}\Big( C_{\boldsymbol a}1_{\{s+L_{\boldsymbol a}\leq t\}}+P_{\boldsymbol a}1_{\{s+L_{\boldsymbol a}> t\}}\Big)y_t+   C_{\boldsymbol a}1_{\{s+L_{\boldsymbol a}\leq T\}}z,
   \end{align*}
where
 \begin{align*}
C_{\boldsymbol a}&=g_0+(T-s-L_{\boldsymbol a}) c+  \sum_{j\in\gamma}g^j +\sum_{j\notin \gamma} (h^j+(L_{\boldsymbol a}+a^j)m^j)\wedge b^j_{(s+L_{\boldsymbol a},a^j+L_{\boldsymbol a})},\\
P_{\boldsymbol a}&=h_0+(T-t)c+\sum_{j=1}^{n}(h^j+(t-s+a^j)m^j)\wedge b^j_{(t,a^j+t-s)}.
\end{align*}

Notice that the total cost function $F_{(s,\boldsymbol a)}(\boldsymbol y_s,z)$ does not explicitly depend on $\boldsymbol x_s$. The role of $\boldsymbol x_s$ becomes explicit through the following additional constraint 
\begin{align}
    (h^j+(a^j+t-s)m^j)x_t^j+b^{j}_{(t,a^j+t-s)}(y_t-x_t^j)=B^j_{(t,a^j+t-s)}y_t,\nonumber\\
    \quad t=s+1,\ldots T,\ j=1,\ldots,n.\label{x_j}
\end{align}
It says that if  $y_t=1$, that is if a PM for at least one component is scheduled at time $t$, then for each component $j$, there is a choice between two actions at time $t$: either perform a PM, so that $x_t^j=1$ and $y_t-x_t^j=0$, or do not perform a PM  and compensate for the current age of the component by increasing the cost function using the virtual replacement cost value (corresponds to $x_t^j=0$ and $y_t-x_t^j=1$).

The optimal maintenance plan is the solution of the linear optimization problem
\begin{align*}
 \text{minimize}\qquad & f_{(s,\boldsymbol a)}(\boldsymbol y_s,z)=\rE(F_{(s,\boldsymbol a)}(\boldsymbol y_s,z))\\
\text{subject to}\qquad & \text{linear constraints } \eqref{tu},\ \eqref{y} \ \eqref{xy} \ \text{ and } \eqref{x_j},\\
&x_t^j\in \{0,1\},\quad t=s+1,\ldots T,\ j=1,\ldots,n,\\
&y_t\in \{0,1\},\quad t=s+1,\ldots T,\\
&z\in \{0,1\}.
\end{align*}

%\section{An opportunistic maintenance plan}
%{\color{red}If any of the components breaks down before the planned next PM, a CM replacement alongside with opportunistic replacements are performed. The opportunistic replacement work as follows: since the maintenance personal need to go there and perform CM on the broken component, they may as well maintain other components if they are close to break down to save the logistic cost. So, for each other component $j$, we compare the virtual maintenance cost $b^{j}_{(t,a^j+t-s)}$ and the PM cost $h^j+(a^j+t-s)m^j$, if virtual maintenance cost is higher, it means that component is too old, it is beneficial to perform PM on the corresponding component. 
%}
\section{Case studies}\label{CS}
This section contains five computational studies with  our model applied to a four component model of a wind power turbine described next. 
Table \ref{tab} lists the four components in question and summarises the basic values of the model parameters, where the suggested values of $(g^j,\beta^j,\theta^j)$  are taken from the paper \cite{tian2011condition}.
\begin{small}
\begin{table}[ht]
\newcommand{\tabincell}[2]{\begin{tabular}{@{}#1@{}}#2\end{tabular}}
  \centering
  \begin{tabular}{|l|c|c|c|c|}\hline
 \small Component ($j$)&\tabincell{c}{\small $g^j$ , \small CM\\ \small replacement\\ \small cost  $[1000 \$]$}&
  \tabincell{c}{\small $m^j$, value  loss \\\small per month \\\small  $[1000 \$]$}&\tabincell{c}{\small Weibull \\\small shape\\\small 
    $\beta^j$}&\tabincell{c}{\small Weibull\\\small scale  
    \\ \small $\theta^j$}\\\hline \small Rotor $(j=1)$&162&0.35&3&1e-06\\\small Main
 Bearing $(j=2)$&110&0.2&2&6.4e-05\\\small Gearbox $(j=3)$&202&0.4&3&1.95e-06\\ \small Generator $(j=4)$&150&0.3&2&8.26e-05\\\hline
\end{tabular}
  \caption{Base values for the $n=4$ component model.}
  \label{tab}
\end{table} 
\end{small}
The cost unit is $1000 \$$, and the time unit is one  month.  In accordance with paper \cite{ZIEGLER20181261}, the lifetime of the wind turbine is assumed to be $20$ years,  so that $T=240$. Other basic values of the model are 
\[g_0=10,\quad h_0=10,\quad h^1=45,\quad h^2=30,\quad h^3=60,\quad h^4=40.\]
%We interpret these cost parameters as follows:
%\begin{description}
%\item[ ] $g_0$ is a sared CM cost,
%\item[ ] $g^j$ is a component specific CM cost,
%\item[ ] $h_0$ is a shared PM cost,
%\item[ ] $h^j+am^j$ is a PM cost assigned to component $j$ having age $a$.
%\end{description}

\subsection{Sensitivity analysis 1}
In this section, we focus on the generator component of the four component model and have a closer look at the cost function of age $a$
$$B_a^4=(h^4+am^4)\wedge b^4_a$$
for different  values of the parameters $(g, h_0,h^4,m^4)$ keeping unchanged the announced values  for the  other parameters of the model .
\begin{figure}[H]
\centering
    \includegraphics[width=12 cm]{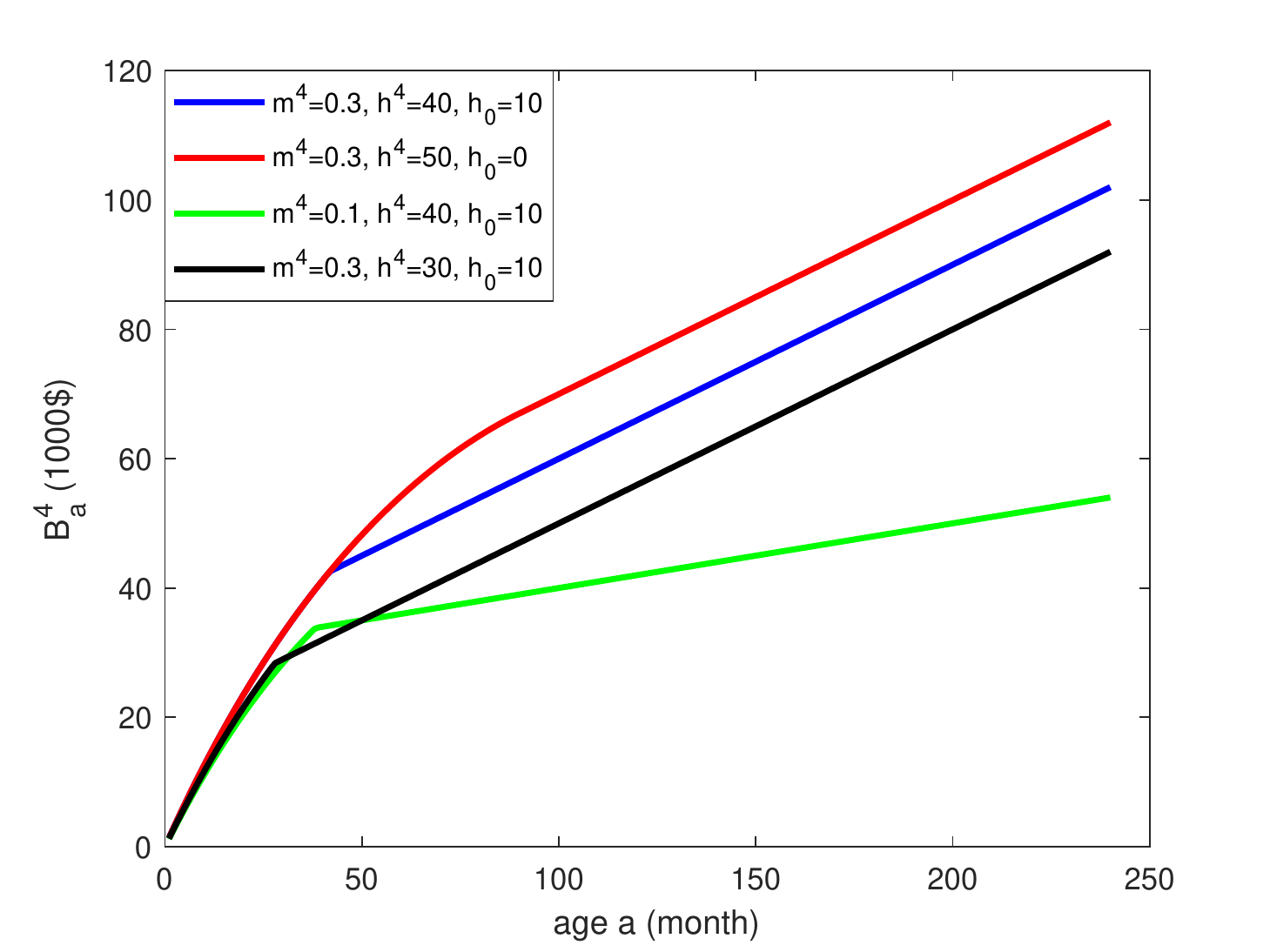}
    \caption{Plots of $B_a^4$ for different combinations of parameters $(g, h_0,h^4,m^4)$.}
    \label{fig1}
\end{figure}

Figure \ref{fig1} summarises the results of this sensitivity analysis. The blue line on Figure \ref{fig1} describes the function $B_a^4$ of age for the base line set of parameter values. For small values of the age variable $a$, the cost function $B_a^4=b^4_a$ grows in a concave manner and beyond some critical age takes the linear form $B_a^4=h^4+am^4$.

The red line shows what happens if we reduce $h_0$ but keep the sum $h^4+h_0$ unchanged. We see that the initial concave part of the curve is the same but the critical age becomes larger. 

The green line indicates a drop in the cost function $B_a^4$ in response to the lowering of the monthly value loss parameter $m^4$. Finally,  the black line shows the cost reduction caused by a smaller value of $h^4$. The black, blue, and red straight lines have the same slope because of the shared parameter value $m^4=0.3$.

\subsection{Sensitivity analysis  2}
One of the most crucial parameters of our model is $m$, the monthly value depreciation of a generic component. In this section, we consider the one component model taking the rotor component as an example. We study how the optimal time to perform the next PM increases as $m=m^1$ becomes larger. All other parameters are fixed at their primary values.
\begin{figure}[H]
\begin{minipage}[t]{0.48\linewidth}
    \includegraphics[width=1\textwidth]{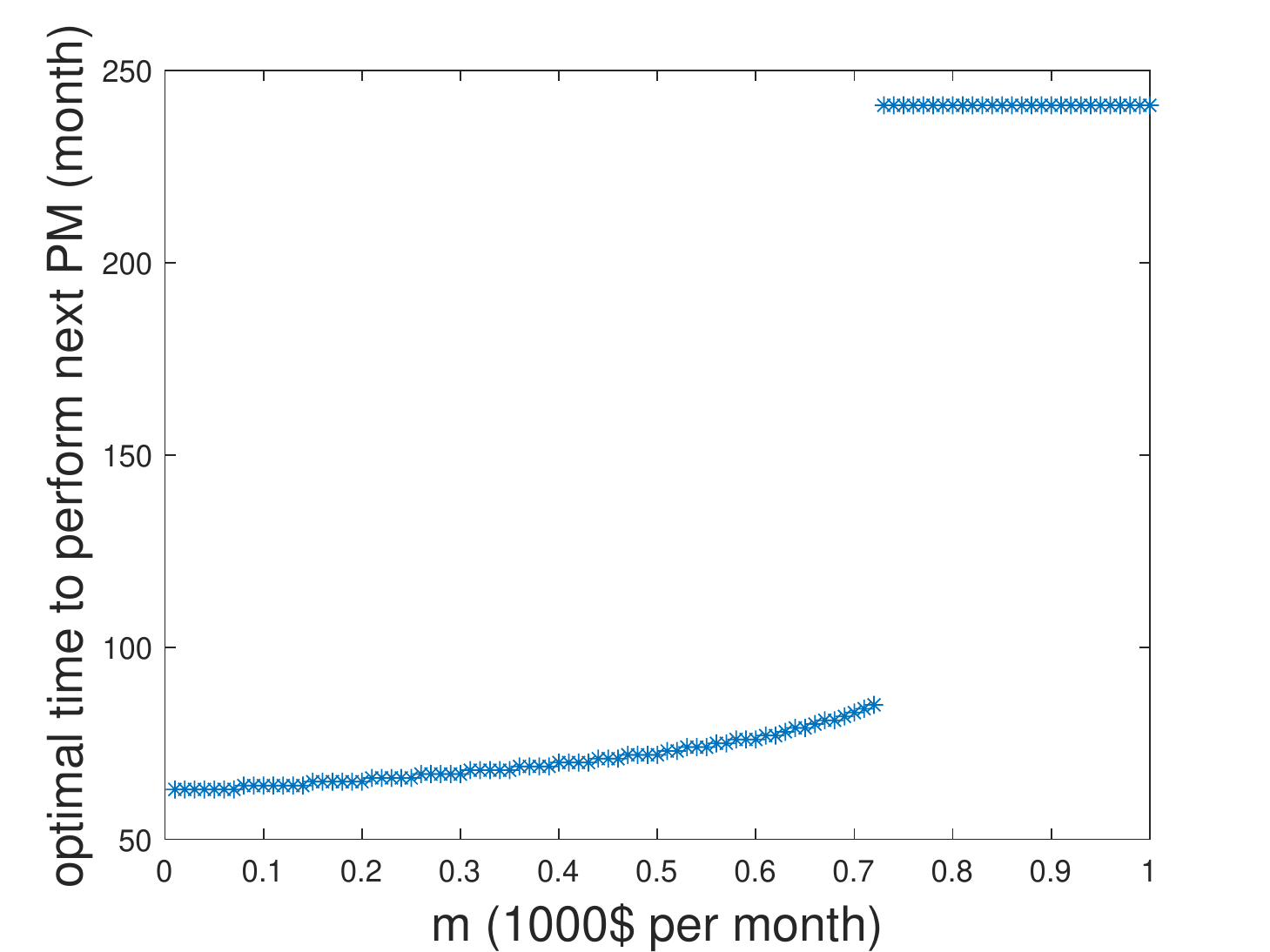}
    \end{minipage}
\begin{minipage}[t]{0.48\linewidth}
\centering
\includegraphics[width=1\textwidth]{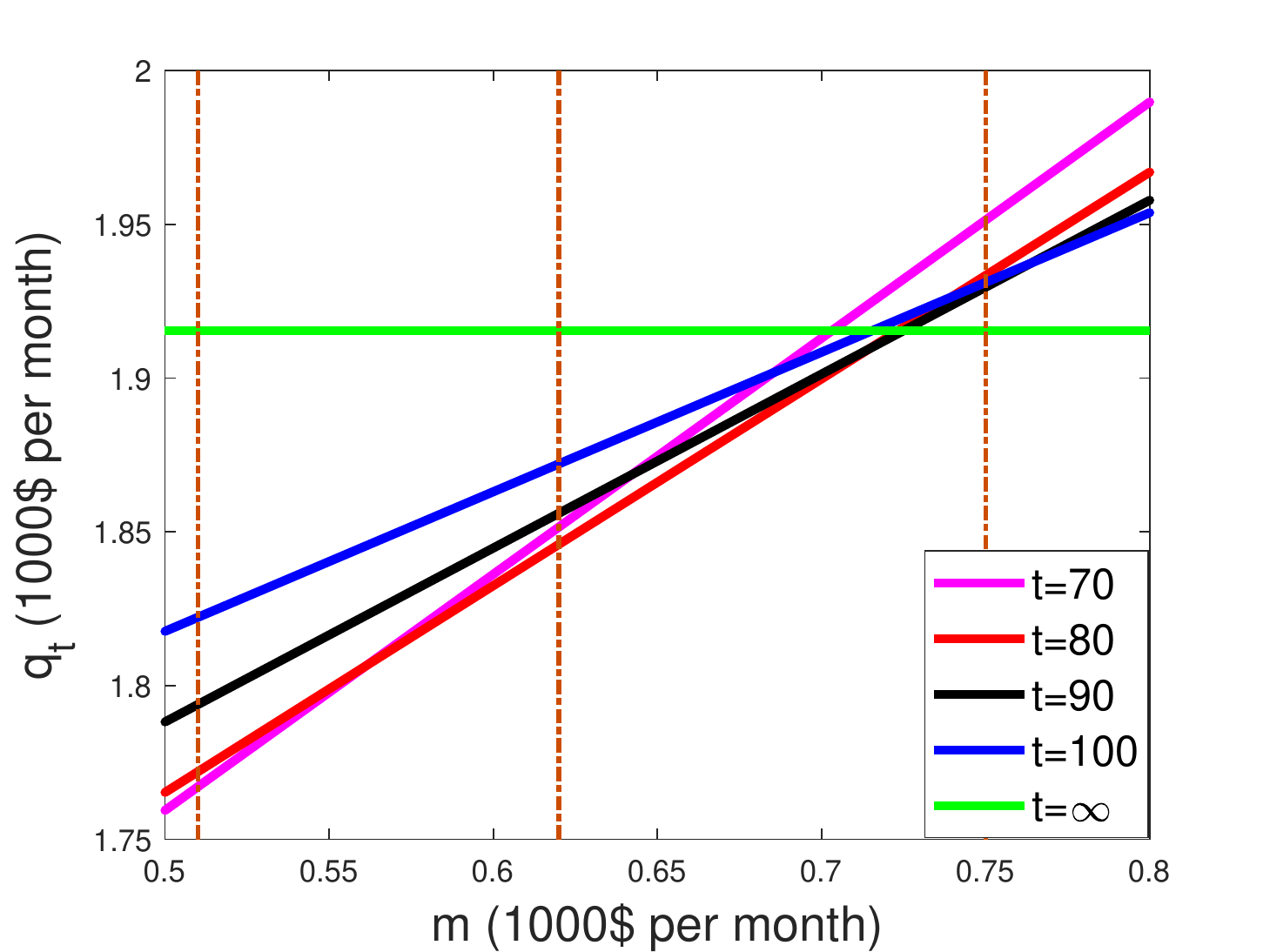}
    \end{minipage}
    \caption{Left panel: the optimal time to perform the next PM as a function of the parameter $m$. Right panel: different average cost based on different PM plan and different monthly value loss $m$.}
    \label{fig2}
\end{figure}

Recall the formula for the time-average maintenance cost in the one component setting: 
\[q_t=\frac{(1-e^{-\theta t^\beta})g+e^{-\theta t^\beta}(h+mt)}{\sum_{k=1}^t (e^{-\theta (k-1)^\beta}-e^{-\theta k^\beta})k+e^{-\theta t^\beta}t},\]
where  $t$ is the planning time for the next PM. Considered as a function of $m$, the average cost $q_t$ is a linear function as depicted on the right panel of Figure \ref{fig2} for different values of the planning time $t$.

The left panel of Figure \ref{fig2} reveals an interesting phenomenon: there exists a critical value $m_0$ of the parameter $m$, beyond which the optimal PM plan is to never perform a PM activity. Indeed, at the value $m_0=0.72$ we observe a jump from optimal planning time $t=80$ up to $t=T$ which is equivalent to $t=\infty$ of no planned PM.

The sudden jump on the left panel graph is explained on the right panel by comparing $q_t$ for different values of $t$. The five straight lines compare $q_{70}, q_{80}, q_{90}, q_{100}, q_{\infty}$ with respect to different values of the parameter $m$. For $m=0.51$, the minimal cost among the five options is given by $q_{70}$. For $m=0.62$, the minimal cost is given by $q_{80}$.
For $m=0.75$, the minimal cost is given by $q_{\infty}$. It is also clear that starting for all $m> 0.73$ the horizontal line, corresponding to $t=\infty$, is always giving the minimal cost.

\subsection{Sensitivity analysis 3}
\begin{figure}[H]
\centering
    \includegraphics[width=6 cm]{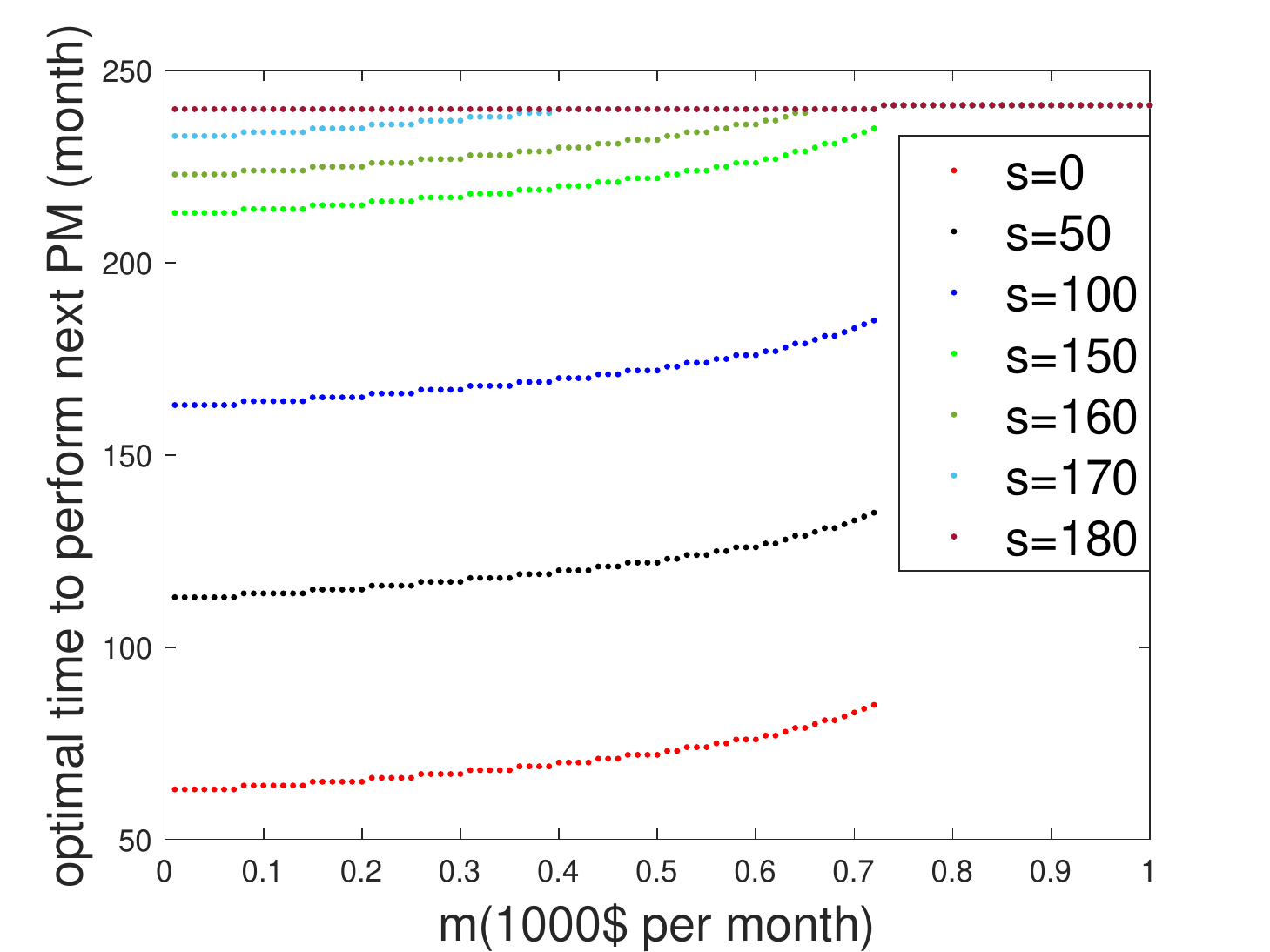}    \includegraphics[width=6 cm]{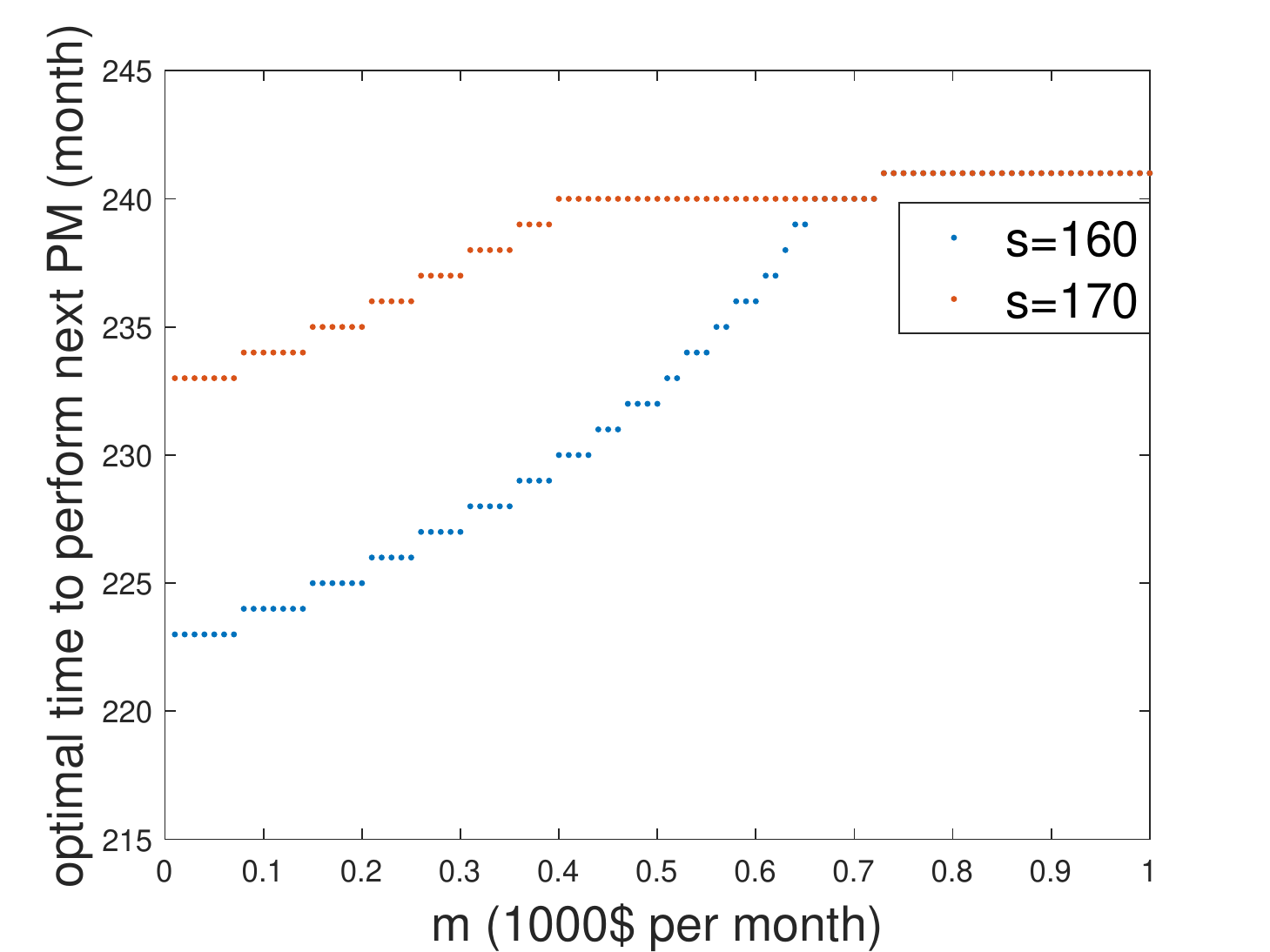}
    \caption{Plots of $t_{(s,0)}$ for different  $s$.}
    \label{fig5}
\end{figure}

In this section, we illustrate the relation \eqref{tas} for the one component (rotor). The Figure \ref{fig5}  compares the optimal times $t_{(s,0)}$ for different starting times $s$ of the planning period.
The curves confirm the stated relation $t_{(s,0)}=t_{(0,0)}+s$. They also emphasise that for $s$ closer to $T$, the approximation  $(T-s)c$ of the maintenance costs  in \eqref{Qsa} becomes rather hoarse for our model to give a reasonable answer. A recent paper \cite{paper2} proposes an optimization model with a special attention to the planing period near to the end time $T$.

However, relation \eqref{tas} allows an effective time effective implementation of our algorithm as a key ingredient of an app (for $s$ not too close to  $T$) by pre-calculating the constants $t_{(0,0)}^j$.

\subsection{Case study 1}

This section deals with the full model of four components. 
The results for different initial ages of four components are shown in Table \ref{4comp}.
\begin{table}[h]
\centering
  \begin{tabular}{|c|c|c|c|c||c|c|}\hline
Initial ages&$j=1$&$j=2$&$j=3$&$j=4$&PM&CM\\\hline 
$(0,0,0,0)$&x&x&61&x&6.06 &6.61 \\
$(30,30,30,30)$ &x&x&31&x&6.59&7.21  \\
$(30,30,0,30)$&39&x&x&x&6.40&6.97 \\
$(40,90,30,60)$ &27&27&27&27&6.71&7.36  \\\hline 
\end{tabular}
\caption{The next PM plan for $4$ components with different initial ages}
\label{4comp}
\end{table}
Comparing the first two cases, we find that $61=30+31$, in accordance with  Proposition 2. 
The two rightmost columns compare monthly maintenance cost for PM planning and pure CM strategy. For the pure CM strategy, we consider the simplest wind turbine maintenance strategy when the PM option is ignored and a CM activity is performed whenever a turbine component breaks down. The pure CM monthly costs are obtained by choosing a parameter $m$ value so high that it is never beneficial to plan a PM activity.
Note that the maintenance cost becomes larger if the components are older at the start. On the other hand, according to Table \ref{4comp} the PM planning may save around 1500 of US dollars per month compared to the pure CM strategy. 

\subsection{Case study 2}
Here the optimization model developed in this paper is compared  with the NextPM model from paper \cite{paper1}. To adapt to the assumption of age independent PM costs of paper \cite{paper1}, we set $m^j=0$ for all $j$. Put $h^1=36.75$, $h^2=23.75$, $h^3=46.75$, $h^4=33.75$, and $g_0=h_0=d$ for the values $d=1,5,10$ considered in paper \cite{paper1}.

The following three tables juxtapose the results produced by the two methods:\\\\
\small 
 \noindent  \begin{tabular}{l|c|c|c|c||c|r}
$d=1$&1&2&3&4&Monthly maintenance cost&Matlab\\\hline 
NextPM&x&x&43&x&4.731 &49 sec\\
New model &x&x&43&x&4.703 &2 sec%\\\hline
\end{tabular}

\vspace{0.5cm}

 \noindent  \begin{tabular}{l|c|c|c|c||c|r} 
$d=5$&1&2&3&4&Monthly maintenance cost&Matlab\\\hline 
NextPM&50&50&50&50&4.964 &54 sec\\
New model&51&51&51&51&4.881  &2 sec%\\\hline
\end{tabular}

\vspace{0.5cm}

\noindent   \begin{tabular}{l|c|c|c|c||c|r}
$d=10$&1&2&3&4&Monthly maintenance cost&Matlab\\\hline 
NextPM&52&52&52&52&5.061 &55 sec\\
New model&52&52&52&52&5.040  &2 sec%\\\hline
\end{tabular}
\normalsize

\vspace{0.5cm}

The optimal schedules are almost identical, demonstrating the accuracy of the new method. Importantly, the new model only requires Matlab to solve it, and the new algorithm is much faster than the previous one.

\section{Conclusions}
In this paper, we studied the scheduling problem for wind turbine maintenance. For maintenance activities, alongside with the traditional CM , PM, and opportunistic replacements, we introduced a new concept of a virtual replacement for taking account of hidden future costs associated with positive ages of components which are treated by the model as good as new. 
We developed an optimization model based on the renewal-reward theorem (which was the main inspiration for the idea of  the virtual replacement costs). 

We tested our optimization model using two case studies in the four-component setting. The first case study showed that with time dependent PM costs, the choice of the initial ages of the components affects the optimal next PM plan significantly. Compared to the pure CM strategy, the first case study produced cost savings close to  $30\%$. In the second case study, we compared the proposed model with an earlier optimization model. The comparison demonstrated that this model is both fast and accurate, so that it can be used as an ingredient of an efficient  maintenance scheduling app for wind power systems.

\section*{Acknowledgements}
\noindent We acknowledge
the financial support from the Swedish Wind Power Technology Centre at Chalmers, the Swedish Energy Agency and Västra Götalandsregionen.

\appendix 
\section{Proofs of propositions}\label{proo}
\subsection*{Proof  of Proposition \ref{propq}}

We will need the following result from the renewal theory, see for example  \cite{grimmett2020probability}. Let 
$$(X_i,R_i),\quad i = 1,2,\ldots$$ 
be independent and identically distributed  pairs of possibly dependent random variables: $X_i$ are positive inter-arrival
times and $R_i$ are associated rewards. Define the cumulative reward process by 
\[W(u) = R_1 + \ldots + R_{N(u)},\]
where $N(u)$ is the number of renewal events up to time $u$, that is $N(u)=k$ if
\[X_1+\ldots+X_k\le u<X_1+\ldots+X_{k+1}.\] According to the renewal-reward theorem, the per unit of time reward 
\[\tfrac{W(u)}{u}\to \tfrac{\rE(R)}{\rE(X)},\quad u\to\infty \]
converges almost surely.

Recall that the full maintenance cost of the gearbox due to a CM is
$g$,
and a similar cost for PM is
$(ma+h)$, where $a$ is the the gearbox age at the replacement. 
Thus, if we plan our next PM for a new gearbox at time $t$ after each renewal event, then the corresponding renewal-reward process has inter-arrival time $X=X(t)$ with
\begin{equation}
X(t)=L\wedge t=L\cdot 1_{\{L\leq t\}}+ t\cdot 1_{\{L\ge t+1\}} 
\label{X0}
\end{equation}
and the reward $R=R(t)$ with
\begin{equation}
R(t)=g 1_{\{L\leq t\}}+ (tm+h) 1_{\{L\ge t+1\}}.
\label{R0}
\end{equation}
Since 
\begin{align*}
 \rE(X)&=\rE(L\cdot 1_{\{L\leq t\}})+ t\rP(L\ge t+1)=\sum_{k=1}^t (e^{-\theta (k-1)^\beta}-e^{-\theta k^\beta})k+e^{-\theta t^\beta}t,\\
 \rE(R)&=g\rP(L\leq t)+ (tm+h)\rP(L\ge t+1)=(1-e^{-\theta t^\beta})g+e^{-\theta t^\beta}(tm+h),
\end{align*}
the renewal-reward theorem implies that the time-average maintenance cost is computed as the following function of the planning time $t$
\[q_t=\frac{\rE(R)}{\rE(X)}=\frac{(1-e^{-\theta t^\beta})g+e^{-\theta t^\beta}(tm+h)}{\sum_{k=1}^t (e^{-\theta (k-1)^\beta}-e^{-\theta k^\beta})k+e^{-\theta t^\beta}t}.\]

\subsection*{Proof  of Proposition \ref{propc}}
By the definition of $t_{(s,a)}$ we have
\[\rE(Q_{(s,a)}(t_{(s,a)},s+L_a))\le \rE(Q_{(s,a)}(t,s+L_a)),\quad t=s+1,\ldots, T+1.\]
To prove the assertion, it is sufficient to verify that
\begin{align}
 \rE(Q_{(s+\delta,a+\delta)}(t_{(s,a)},s+\delta+L_{(a+\delta)}))&\le \rE(Q_{(s+\delta,a+\delta)}(t,s+\delta+L_{(a+\delta)})), \label{ineq}
 \end{align}
for $t=s+\delta+1,\ldots, T+1$.

To derive \eqref{ineq}, we will introduce special notation
\begin{align*}
 B(t)&=g+(T-t)c,\\
\hat  B(t,u)&=tm+h+(T-u)c,
\end{align*}
and notice that for $t\le T$
\begin{align*}
\rE(Q_{(s,a)}(t,s+L_a))&= \sum_{u=1}^{t-s}B(s+u)\rP(L-a=u|L> a)\\
&+ \hat B(t-s+a,t)\rP(L-a> t-s|L>a)
\end{align*}
or in other words,
\begin{align*}
\rP(L> a)\rE(Q_{(s,a)}(t,s+L_a))&= \sum_{u=1}^{t-s}B(s+u)\rP(L=a+u)\\
&+ \hat B(t-s+a,t)\rP(L> a+t-s).
\end{align*}
On the other hand, we have similarly,
\begin{align*}
\rP&(L> a+\delta)\rE(Q_{(s+\delta,a+\delta)}(t,s+\delta+L_{(a+\delta)}))\\
&= \sum_{u=1}^{t-s-\delta}B(s+\delta+u)\rP(L=a+\delta+u)+ \hat B(t-s+a,t)\rP(L> a+t-s).
\end{align*}
The key observation leading to \eqref{ineq} is that the difference
 \begin{align*}
\rP(L> a)\rE(Q_{(s,a)}(t,s+L_a))&-\rP(L> a+\delta)\rE(Q_{(s+\delta,a+\delta)}(t,s+\delta+L_{(a+\delta)}))\\
&= \sum_{u=1}^{\delta}B(s+u)\rP(L=a+u)
\end{align*}
is independent of $t$. In view of this fact, it is clear that the earlier observed inequality
\[\rP(L> a)\rE(Q_{(s,a)}(t_{(s,a)}, s+L_a))\le \rP(L> a)\rE(Q_{(s,a)}(t,s+L_a))\]
entails
\begin{align*}
  \rP(L> a+\delta)\rE(&Q_{(s+\delta,a+\delta)}(t_{(s,a)},s+\delta+L_{(a+\delta)}))\\
  &\le  \rP(L> a+\delta)\rE(Q_{(s+\delta,a+\delta)}(t,s+\delta+L_{(a+\delta)})),
\end{align*}
which immediately implies \eqref{ineq}.

\subsection*{Proof  of Proposition \ref{prop6}}

If $t\le T$, then according to \eqref{Qsa}
\[
        Q_{(s,0)}(t,L)=[g+(T-s-L)c]\cdot 1_{\{s+L\leq t\}}  +  [h+(t-s)m+(T-t)c]\cdot 1_{\{s+L>t\}},
\]
and 
\begin{align*}
        \rE(Q_{(s,0)}(t,L))-(T&-s)c=g\rP(L\leq t-s)+(h+(t-s)m)\rP(L>t-s)\\
        &\hspace{2cm}-c(\rE(L1_{\{L\leq t-s\}})+(t-s)\rP(L>t-s))\\
        &=\rE(R(t-s))-c\rE(X(t-s))=(q_{t-s}-c)\rE(X(t-s)),
\end{align*}
where in the last line we use notation \eqref{X0} and \eqref{R0}.
Since by definition, $c\le q_{t-s}$, the obtained equality
\[\rE(Q_{(s,0)}(t,L))=(T-s)c+(q_{t-s}-c)\rE(X(t-s))\]
implies that provided $t_{(0,0)}+s\le T,$
$$f^*_{(s,0)}=\min_t\rE(Q_{(s,0)}(t,L))=(T-s)c,$$ 
and we conclude 
\begin{equation}
f^*_{(s_1,0)}-f^*_{(s_2,0)}=(s_2-s_1)c.
\label{bb}
\end{equation}
Furthermore, we see that
\[t_{(s,0)}=t_{(0,0)}+s,\]
which together with Proposition \ref{propc} yields \eqref{tas}.

It remains to prove that $b_{(s,a)}=b_a$ or equivalently,
\[f_{(s,a)}^*=f_{(0,a)}^*-sc.\]
Since
\begin{align*}
f_{(s,a)}^*&=\rE(Q_{(s,a)}(t_{(s,a)},s+L_a))=\rE(Q_{(s,a)}(t_{(0,a)}+s,s+L_a))\\
&=\rE\Big((g+(T-s-L_a)c)1_{\{L_a\leq t_{(0,a)}\}}\\
&\hspace{2cm} +(h+m(t_{(0,a)}+a)+(T-t_{(0,a)}-s)c)1_{\{L_a>t_{(0,a)}\}}\Big)\\
&=\rE(Q_{(0,a)}(t_{(0,a)},L_a))-sc=f_{(0,a)}^*-sc.
\end{align*}
This finishes the proof of Proposition \ref{prop6}.

%\bibliographystyle{elsarticle-num} 

\bibliographystyle{plain}
\bibliography{elsarticle-template}
\end{document}